\input amstex
\documentstyle {amsppt}
\pagewidth{12.5 cm}\pageheight{19 cm}\magnification\magstep1
\topmatter   
\title Fermionic form and Betti numbers\endtitle
\author G. Lusztig\endauthor
\abstract{We state a conjectural relationship between the fermionic form 
\cite{HKOTY} and the Betti numbers of a Grassmannian over a preprojective
algebra or, equivalently, of a lagrangian quiver variety.}\endabstract
\address Department of Mathematics, M.I.T., Cambridge, MA 02139\endaddress
\thanks{Supported in part by the National Science Foundation}\endthanks
\endtopmatter
\document

\define\mto{\mapsto}

\define\tim{\times}

\redefine\i{^{-1}}
\define\fra{\frac}

\define\ot{\otimes}

\define\bcc{\thickfracwithdelims[]\thickness0}

\define\opl{\oplus}

\define\al{\alpha}

\define\de{\delta}

\define\om{\omega}

\define\la{\lambda}

\define\vp{\varpi}
\define\vt{\vartheta}

\define\bc{\bold C}
\define\bd{\bold D}

\define\bi{\bold I}

\define\bn{\bold N}

\define\bp{\bold P}

\define\br{\bold R}

\define\bz{\bold Z}

\define\cf{\Cal F}

\define\ci{\Cal I}

\define\cv{\Cal V}

\define\GP{GP}
\define\HKOTY{HKOTY}
\define\KR{KR}
\define\KI{Ki}
\define\KSS{KSS}
\define\KLE{Kl}
\define\LU{L1}
\define\LUU{L2}
\define\LUUU{L3}
\define\NA{NA}
\define\Gra{\text{\rm Grass}}
\subhead 1. Notation \endsubhead
We fix a graph of type $ADE$ with set of vertices $I$. Let $E$ be a 
$\br$-vector space with a basis $(\al_i)_{i\in I}$ and a positive definite
symmetric bilinear form $(,):E\tim E@>>>\br$ given by $(\al_i,\al_i)=2$,
$(\al_i,\al_j)=-1$ if $i,j$ are joined in the graph, $(\al_i,\al_j)=0$ if 
$i\ne j$ are not joined in the graph. Let $(\vp_i)_{i\in I}$ be the basis of 
$E$ defined by $(\vp_i,\al_j)=\de_{i,j}$. For $\xi\in E$ define 
${}^i\xi,{}_i\xi$ in $\br$ by
$$\xi=\sum_i({}^i\xi)\vp_i=\sum_i({}_i\xi)\al_i.$$
Let $P=\{\xi\in E|{}^i\xi\in\bz\quad\forall i\in I\}$,
$P^+=\{\xi\in E|{}^i\xi\in\bn\quad\forall i\in I\}$. Let 
$\rho=\sum_i\vp_i\in P^+$. We consider the usual partial order on $P$:
$$\xi\le\xi' \Leftrightarrow {}_i\xi'-{}_i\xi\in\bn \text{ for all } i.$$
For $i\in I$ define $s_i:E@>>>E$ by $s_i(\xi)=\xi-(\xi,\al_i)\al_i$. Let $W$ be
the (finite) subgroup of $GL(E)$ generated by $\{s_i|i\in I\}$. Let $\bz[P]$ be
the group ring of $P$ with obvious basis $([\xi])_{\xi\in P}$. For $\xi\in P^+$
define $V_\xi\in\bz[P]$ by Weyl's formula
$$\sum_{w\in W}\det(w)[w(\xi+\rho)]=V_\xi\sum_{w\in W}\det(w)[w(\rho)].$$

\subhead 2. The fermionic form \cite{\HKOTY}\endsubhead
Let $q$ be an indeterminate. For $p,m\in\bn$ define 
$$\bcc{p+m}{m}=\fra{(q^{p+1}-1)(q^{p+2}-1)\dots(q^{p+m}-1)}
{(q-1)(q^2-1)\dots(q^m-1)}\in\bz[q].$$
Let $\nu=\{\nu_k^{(i)}\in\bn|i\in I, k\ge 1)$ where all but finitely many 
$\nu_k^{(i)}$ are zero. Let $\la\in P^+$. In \cite{\HKOTY, 4.3} a "fermionic 
form" $M(\nu,\la,q)$ (or $M(W,\la,q)$ in the notation of {\it loc.cit.}) is
attached to $\nu,\la$. This is a $q$-analogue of an expression which first 
appeared in Kirillov and Reshetikhin \cite{\KR}. For $q=1$ it conjecturally 
gives the multiplicities in certain representations of an affine quantum group
when restricted to the ordinary quantum group. In \cite{\KLE}, Kleber rewrites
the formula of \cite{\KR} in the form of a computationally efficient algorithm.
(In his paper, it is assumed that one of the $\nu_k^{(i)}$ is $1$ and the 
other s are $0$ but, as he pointed out to me, the same procedure works in 
general for $q=1$.)

In the remainder of this note we assume that
$$\nu_k^{(i)}=0 \text{ for } i\in I, k\ge 2.$$
In this case we identify $\nu$ with the element of $P^+$ such that 
${}^i\nu=\nu_1^{(i)}$ for all $i$. By definition,
$$M(\nu,\la,q)=\sum_{\{m\}}q^{c(\{m\})}\prod_{i\in I;k\ge 1}
\bcc{p_k^{(i)}+m_k^{(i)}}{m_k^{(i)}},\tag a$$
$$c(\{m\})=\fra{1}{2}\sum_{i,j\in I}(\al_i,\al_j)\sum_{k,l\ge 1}
\min(k,l)m_k^{(i)}m_l^{(j)}-\sum_{i\in I}\sum_{k\ge 1}{}^i\nu m_k^{(i)},$$
$$p_k^{(i)}={}^i\nu-\sum_{j\in I}(\al_i,\al_j)\sum_{l\ge 1}\min(k,l)m_l^{(j)},
$$
where the sum $\sum_{\{m\}}$ is taken over $\{m_k^{(i)}\in\bn|i\in I, k\ge 1\}$
satisfying $p_k^{(i)}\ge 0$ for $i\in I, k\ge 1$ and
$$\sum_{i\in I}\sum_{k\ge 1}km_k^{(i)}\al_i=\nu-\la$$
for $i\in I$. We rewrite this by extending the method of \cite{\KLE} to the 
$q$-analogue; we obtain
$$M(\nu,\la,q)=\sum_\om q^{c(\om)}\prod_{i\in I;k\ge 1}
\bcc{{}^i\om_k+{}_i\mu_k}{{}_i\mu_k}\tag b$$
sum over all sequences $\om$ in $P^+$ of the form
$$\nu=\om_0>\om_1>\om_2>\dots>\om_s=\om_{s+1}=\om_{s+2}=\dots=\la$$
such that 
$$\om_0-\om_1\ge\om_1-\om_2\ge\om_2-\om_3\ge\dots$$ 
that is,
$$\mu_k=\om_{k-1}-2\om_k+\om_{k+1}\ge 0 \text{ for } k\ge 1,\quad \mu_k=0
\text{ for } k\gg 0,$$ 
and
$$c(\om)=\fra{1}{2}\sum_{k\ge 1}(X_k,X_k)-(\nu,X_1)$$
where
$$X_k=\om_{k-1}-\om_k \text{ for } k\ge 1.$$
The connection between (a) and (b) is as follows: in terms of the data in (a) 
we have
$$\om_k=\sum_ip_k^{(i)}\vp_i,\quad \mu_k=\sum_im_k^{(i)}\al_i.$$
Since $\mu_k=X_k-X_{k+1}$ for $k\ge 1$, we have for $i,j\in I$:
$$\align&\sum_{k,l\ge 1}\min(k,l)m_k^{(i)}m_l^{(j)}=
\sum_{k,l\ge 1}\min(k,l)({}_iX_k-{}_iX_{k+1})({}_jX_l-{}_jX_{l+1})\\&
=\sum_{k,l\ge 1}\min(k,l)({}_iX_k({}_jX_l)-{}_iX_{k+1}({}_jX_l)
-{}_iX_k({}_jX_{l+1})+{}_iX_{k+1}({}_jX_{l+1}))\\&
=\sum_{k,l\ge 1}(\min(k,l)-\min(k-1,l)-\min(k,l-1)+\min(k-1,l-1))
{}_iX_k({}_jX_l)\\&=\sum_{k\ge 1}{}_iX_k({}_jX_k),\endalign$$
$$\sum_{k\ge 1}{}^i\nu m_k^{(i)}=\sum_{k\ge 1}{}^i\nu({}_iX_k-{}_iX_{k+1})
=({}^i\nu)({}_iX_1),$$
hence $c(\{m\})=c(\om)$.

The following result is stated without proof in \cite{\HKOTY}.
\proclaim{Lemma 3} $M(\nu,\la,q)\in\bn[q\i]$.
\endproclaim
Let $\om$ be as in Sec.2. The product of $q$-binomial coefficients in the term
corresponding to $\om$ is a polynomial in $q$ of degree
$$\align&N=\sum_{i\in I;k\ge 1}{}^i\om_k({}_i\mu_k)=
\sum_{k\ge 1}(\om_k,\mu_k)\\&=\sum_{k\ge 1}(\nu-X_1-X_2-\dots-X_k,X_k-X_{k+1})
=(\nu,X_1)-\sum_{k\ge 1}(X_k,X_k).\endalign$$
It is enough to show that $c(\om)+N\le 0$. We have
$$c(\om)+N=-\fra{1}{2}\sum_{k\ge 1}(X_k,X_k)$$ 
and this is clearly $\le 0$.

\proclaim{Lemma 4}Let $\xi\in P^+$ and let $\eta\in P$ be such that 
$\eta\ge 0$. Then $(\xi,\eta)\ge 0$.
\endproclaim
This is obvious.

\proclaim{Lemma 5} If $\nu\ge\la$ then
$M(\nu,\la,q)=q^{-(\nu,\nu)/2+(\la,\la)/2}+$ strictly larger powers of $q$.
\endproclaim
Let $\om$ be as in Sec.2. We show that
$$c(\om)\ge-(\nu,\nu)/2+(\la,\la)/2 \tag a$$
that is,
$$\fra{1}{2}\sum_{k\ge 1}(X_k,X_k)-(\nu,X_1)-\fra{1}{2}(\nu-\la,\nu-\la)+
(\nu,\nu-\la)\ge 0.$$
Since $\nu-\la=X_1+X_2+X_3+\dots$, this is equivalent to
$$(\nu,X_2+X_3+\dots)-\sum_{1\le k<l}(X_k,X_l)\ge 0.\tag b$$
Applying Lemma 4 to $\xi=\nu-X_1-X_2-\dots-X_{k+1}=\om_{k+1}$, 
$\eta=X_{k+1}$, we obtain
$$(\nu-X_1-X_2-\dots-X_{k+1},X_{k+1})\ge 0.$$
Adding these inequalities over all $k\ge 1$ we obtain
$$\sum_{k\ge 1}(\nu-X_1-X_2-\dots-X_k-X_{k+1},X_{k+1})\ge 0,$$
that is,
$$(\nu,X_2+X_3+\dots)-\sum_{1\le k<l\ge 1}(X_k,X_l)\ge\sum_{k\ge 2}(X_k,X_k)
\ge 0.$$
Thus, (b) hence (a) are proved. This proof shows also that the inequality (a) 
is strict unless $\om$ satisfies $X_2=X_3=\dots=0$. If this last condition is 
satisfied then $\om$ is the sequence $\nu=\om^0\ge\om^1=\om^2=\dots=\la$ and 
(a) is an equality. The lemma is proved.

\subhead 6. Inversion\endsubhead
Define $M^*(\nu,\la,q)\in\bz[q\i]$ for any $\nu,\la\in P^+$ by the requirement
that the matrix $(M^*(\nu,\la,q))_{\mu,\la}$ is inverse to the matrix 
$(M(\nu,\la,q))_{\mu,\la}$ (which is lower triangular with $1$ on diagonal). 
Thus, $M^*(\nu,\nu,q)=1$ and $\sum_{\la\in P^+}M^*(\nu,\la,q)M(\la,\xi,q)=0$ 
for any $\nu>\xi$ in $P^+$. There is some evidence that the matrix $M^*$ is 
simpler than $M$. For example, in type $A_1$, we have
$$M^*(\nu,\la,1)=
(-1)^{{}_i(\nu-\la)}\binom{{}^i\la+{}_i(\nu-\la)}{{}_i(\nu-\la)}$$
for any $\nu\ge\la$ in $P^+$.

\subhead 7. Path algebra \endsubhead
Let $\bi$ be the set of all sequences $i_1,i_2,\dots,i_s$ (with $s\ge 1$) in 
$I$ such that $i_k,i_{k+1}$ are joined for any $k\in[1,s-1]$. Let $\cf$ be the 
$\bc$-vector space spanned by elements $[i_1,i_2,\dots,i_s]$ corresponding to 
the various elements of $\bi$. We regard $\cf$ as an algebra in which the 
product $[i_1,i_2,\dots,i_s][j_1,j_2,\dots,j_{s'}]$ is equal to
$[i_1,i_2,\dots,i_s,j_2,\dots,j_{s'}]$ if $i_s=j_1$ and is zero, otherwise. For
$i\in I$, let $\vt_i=\sum_j [iji]$ where $j$ runs over the elements of $I$ that
are joined with $i$. For $i,j\in I$ let $\cf_{ij}$ be the subspace of $\cf$ 
spanned by the elements $[i_1,i_2,\dots,i_s]$ with $i_1=i,i_s=j$. For $u\in\bz$
let $\cf^u$ be the subspace of $\cf$ spanned by the elements 
$[i_1,i_2,\dots,i_s]$ with $s=u+1$. (For $u<0$ we have $\cf^u=0$.) Let 
$\cf_{ij}^u=\cf_{ij}\cap\cf^u$. We have
$$\cf=\opl_{i,j}\cf_{ij},\cf=\opl_u\cf^u, \cf=\opl_{i,j;u}\cf_{ij}^u.$$
Let $\ci$ be the two-sided ideal of $\cf$ generated by the elements $\vt_i$ 
($i\in I$). The quotient algebra $\bp=\cf/\ci$ has finite dimension over $\bc$
\cite{\GP}. Let $\bp_{ij},\bp^u,\bp_{ij}^u$ be the image of 
$\cf_{ij},\cf^u,\cf_{ij}^u$ in $\bp$. We have 
$$\bp=\opl_{i,j}\bp_{ij}, \bp=\opl_u\bp^u, \bp=\opl_{i,j;u}\bp_{ij}^u.$$
Let $\bd$ a finite dimensional $\bc$-vector with a given direct sum
decomposition $\bd=\opl_{i\in I}\bd_i$. Then  
$\bd^\dag=\opl_{i,j}\bp_{ij}\ot\bd_j$ is a left $\bp$-module in an obvious way
(in fact a projective $\bp$-module of finite dimension over $\bc$). Let 
$\nu=\sum_{i\in I}\dim\bd_i\vp_i\in P^+$. Let $\Gra_\bp(\bd^\dag)$ be the 
algebraic variety consisting of all $\bp$-submodules of $\bd^\dag$. We have a 
partition 
$$\Gra_\bp(\bd^\dag)=\sqcup_{\xi\in P}\Gra_{\bp,\xi}(\bd^\dag)$$
where $\Gra_{\bp,\xi}(\bd^\dag)$ consists of all $\bp$-submodules $\cv$ such
that $\sum_i\dim([i]\bd^\dag)/[i]\cv)\al_i=\nu-\xi$. Then 

\proclaim{Conjecture A} Let $q^{1/2}$ be an indeterminate. For any $\xi\in P$ 
we have
$$\sum_{s\in\bn}\dim H^s(\Gra_{\bp,\xi}(\bd^\dag))q^{s/2}
=\sum_{\la\in P^+}(\xi:V_\la)q^{(\nu,\nu)/2-(\xi,\xi)/2}M(\nu,\la,q)\tag a$$
where $(\xi:V_\la)$ is the coefficient in $\xi$ in $V_\la$ and $H^s()$ denotes
ordinary cohomology with coefficients in a field.
\endproclaim
Since $(\xi:V_\la)$ and $(\xi,\xi)$ are $W$-invariant in $\xi$, we see
that the right hand side of (a) is $W$-invariant in $\xi$. The analogous 
property of the left hand side of (a) is known to be true. (See \cite{\LUU}.)

In \cite{\LU} it is shown that $\Gra_{\bp,\xi}(\bd^\dag)$ is isomorphic to a
(lagrangian) quiver variety defined in Nakajima \cite{\NA} and, conversely, all
such quiver varieties are obtained. Thus the conjecture above gives at the same
time the Betti numbers of quiver varieties.

Assuming the conjecture, we show that $\Gra_{\bp,\xi}(\bd^\dag)$ is connected 
if $(\xi:V_\nu)\ne 0$ (an expected but not yet proved property of quiver 
varieties). We may assume that $\xi\in P^+,\xi\le\nu$. It suffices to show that
$\dim H^0(\Gra_{\bp,\xi}(\bd^\dag))=1$ or that the constant term of
$$\sum_{\la\in P^+}(\xi:V_\la)q^{(\nu,\nu)/2-(\xi,\xi)/2}M(\nu,\la,q)$$
is $1$. By Lemma 5, the constant term of the term corresponding to $\la=\xi$ is
$1$. Consider now the term corresponding to $\la\ne\xi$; we show that its
constant term is $0$. We may assume that $\la\le\nu$ and $(\xi:V_\la)\ne 0$ so
that $\xi<\la$. By Lemma 5, $q^{(\nu,\nu)/2-(\xi,\xi)/2}M(\nu,\la,q)$ is of the
form
$$q^{(\nu,\nu)/2-(\xi,\xi)/2)}q^{-(\nu,\nu)/2+(\la,\la)/2}+ \text{ strictly 
larger powers of } q.$$
Since $(\xi,\xi)<(\la,\la)$ for any $\la,\xi$ in $P^+$ such that $\xi<\la$, our
assertion is established.

The same argument shows that the right hand side of (a) is in $\bn[q]$. 

Assuming again that $\xi\in P^+,\xi\le\nu$ we show that the polynomial in $q$ 
given by the right hand side of (a) has degree $s=(\nu,\nu)/2-(\xi,\xi)/2$. The
term corresponding to $\la=\nu$ is $(\xi:V_\nu)q^s$ where $(\xi:V_\nu)>0$. 
Consider now the term corresponding to $\la\ne\nu$. We may assume that 
$\la<\nu$. It suffices to show that in this case $M(\nu,\la,q)\in q\i\bz[q\i]$.
This follows from the proof of Lemma 3. Our assertion is established. Note that
this is compatible with the conjecture since the dimension of 
$\Gra_{\bp,\xi}(\bd^\dag)$ is known to be equal to $(\nu,\nu)/2-(\xi,\xi)/2$. 

In the $A_1$ case, the left hand side of (a) is a $q$-binomial coefficient; 
the right hand side can be computed by results in \cite{\KI},\cite{\KSS}; the 
conjecture holds in this case.

\subhead 8 \endsubhead
The conjecture implies that 
$$\chi(\Gra_\bp(\bd^\dag))=\sum_{\la\in P^+}\dim(V_\la)M(\nu,\la,1)\tag a$$
where $\chi$ denotes Euler characteristic and 
$\dim(V_\la)=\sum_{\xi\in P}(\xi:V_\la)$. Let $f(\nu)$ (resp. $g(\nu)$) be the
left (resp. right) hand side of (a). According to \cite{\HKOTY}, it is expected
that $g(\nu+\nu')=g(\nu)g(\nu')$ for any $\nu,\nu'\in P^+$. The corresponding 
identity $f(\nu+\nu')=f(\nu)f(\nu')$ is known \cite{\LUUU, 3.20}.

\subhead 9 \endsubhead
We can partition $I$ into two disjoint subsets $I^0,I^1$ so that no two 
vertices in $I^0$ are joined and no two vertices in $I^1$ are joined. For any
$u\in\bz$ let
$${}^u\bd^\dag=\opl_{i\in I^\de,j\in I}(\bp_{ij}^u\opl\bp_{ij}^{u-1})\ot\bd_j$$
where $\de\in\{0,1\}$ is defined by $u=\de\mod 2$. We have 
$\bd^\dag=\opl_u{}^u\bd^\dag$. Consider the $\bc^*$-action $t,d\mto td$ on 
$\bd^\dag$ with weight $u$ on ${}^u\bd^\dag$. This action is compatible with 
the $\bp$-module structure in the following sense: $t(pd)=(tp)(td)$ for 
$t\in\bc^*,p\in\bp,d\in\bd^\dag$ where $tp$ is given by the $\bc^*$-action on 
$\bp$ (through algebra automorphisms) for which $\bp^u$ has weight $u$. Hence
we have an induced $\bc^*$-action on $\Gra_{\bp,\xi}(\bd^\dag)$ for any $\xi$.
Let $\Gra'_{\bp,\xi}(\bd^\dag)$ be the fixed point set of this $\bc^*$-action
(a smooth variety). It consists of all $\cv\in\Gra_{\bp,\xi}(\bd^\dag)$ such
that $\cv=\opl_u(\cv\cap {}^u\bd^\dag)$. 

\proclaim{Conjecture B} For any $\xi\in P$ we have
$$\sum_{s\in\bn}\dim H^s(\Gra'_{\bp,\xi}(\bd^\dag))q^{s/2}
=\sum_{\la\in P^+}(\xi:V_\la)M(\nu,\la,q\i).$$
\endproclaim
One can show that this is equivalent to Conjecture A.

\medpagebreak

I wish to thank M. Kleber and A. Schilling for useful discussions.

\enddocument